\documentclass{article}
\usepackage{authblk}
\usepackage{amsmath,amsthm,amssymb,xcolor,url,booktabs,mathtools}
\usepackage{enumitem,etoolbox}
\usepackage{hyperref}
\usepackage[T1]{fontenc}
\usepackage[textwidth=17cm, textheight=22cm, bottom=30mm]{geometry}

\newtheorem{theorem}{Theorem}
\newtheorem{lemma}[theorem]{Lemma}

\newtheorem{definition}{Definition}
\newtheorem{problem}{Problem}

\newtheorem{remark}[theorem]{Remark}

\newcommand{\BH}{\mathrm{BH}}
\newcommand{\Z}{\mathbb{Z}}

\newcommand{\C}{\mathbb{C}}
\newcommand{\R}{\mathbb{R}}
\newcommand{\F}{\mathbb{F}}
\newcommand{\PAF}{P}
\newcommand{\NAF}{N}

\newcommand{\T}{^{\mathsf T}}

\AtBeginDocument{}

\title{Recent progress in the construction of Hadamard matrices}

\author[1]{Ronan Egan}
\author[2]{Ilias S. Kotsireas}
\author[3]{Stephen London}
\author[4]{P\'adraig \'O Cath\'ain}
\affil[1]{School of Mathematical Sciences, Dublin City University, Ireland}
\affil[2]{Wilfrid Laurier University, Waterloo, Ontario, Canada}
\affil[3]{Loyola University, Chicago, USA}
\affil[4]{\'Udar\'as na Gaeltachta, Ireland}
\date{}

\begin{document}
\maketitle
\textit{To Jennifer R. Seberry and Dragomir \v{Z}. \DJ okovi\'c, pioneers in the existence of Hadamard matrices.}

\begin{abstract}
A Hadamard matrix is a square matrix with entries $\pm 1$ whose rows are mutually orthogonal.
Hadamard matrices are conjectured to exist at every order divisible by $4$; for approximately the last twenty years, the smallest open case of this conjecture was order 668.
In August 2026, Alp\"oge, Voinov and Reynolds-Haertle announced
Hadamard matrices at the twelve orders below $2000$ at which none was known.
In this paper we verify and analyse the claimed matrices. We survey results on the existence of Hadamard matrices and list the $95$ orders up to
$10\,000$ at which no Hadamard matrix is currently known.
In pursuing this work we constructed Hadamard matrices at sixteen orders that were open in the published literature. We conclude with a list of open problems.

MSC: 05B20, 15B34
\end{abstract}

\section{Introduction}

Let $M$ be an $n \times n$ complex matrix with entries of modulus at
most $1$. Hadamard proved that $|\det M| \le n^{n/2}$, with equality if
and only if every entry has modulus $1$ and $MM^{\ast} = nI$
\cite{Hadamard1893}. Over $\C$ the bound is attained at every order by
the character table of $\Z_n$. Over $\R$ the extremal matrices are
Hadamard matrices, which exist only at orders $1$, $2$ and $4t$. The
first infinite family is due to Sylvester \cite{Sylvester1867}. Paley described further infinite families and conjectured that a Hadamard matrix exists at every order $4t$
\cite{Paley1933}.

The conjecture is settled at orders below $92$ by the constructions of
Sylvester and Paley together with Kronecker products. Order $92$ was resolved computationally by Baumert, Golomb and Hall \cite{BaumertGolombHall1962}, and order $116$ by Baumert in 1966 \cite{Baumert116}. By 1977 the smallest open order was $268$ \cite{HedayatWallis1978, Seberry1978}; it was settled by Sawade in 1985 \cite{Sawade268}, leaving $428$ as the smallest unsettled order. From the construction of order $428$ by Kharaghani and Tayfeh-Rezaie \cite{Kharaghani428}, the smallest open order was $668$. On 12 August 2026 Alp\"oge, Voinov and Reynolds-Haertle announced Hadamard matrices at the twelve orders
\[
  668,\ 716,\ 892,\ 1132,\ 1244,\ 1388,\ 1436,\ 1676,\ 1772,\ 1916,\
  1948,\ 1964,
\]
which were the only open orders below $2000$ \cite{Alpoge2026}. The
announcement gives the matrices but not the method of construction.

Verifying the twelve matrices is a matrix multiplication; their structure
is less apparent, and we describe it in detail. Each is a Goethals--Seidel
array over an abelian group. Nine are over cyclic groups, four of them
with a border of width $4$, and three are over non-cyclic groups, with a
wider border that is constant on the cosets of a subgroup.

A Hadamard matrix is therefore known at every order below $2000$ divisible
by $4$. We construct Hadamard matrices at sixteen further orders below
$10\,000$, record the orders settled by others while this paper was in
preparation, and list the $95$ orders $n \le 10\,000$ at which no Hadamard
matrix is currently known; the smallest is $2396 = 4 \cdot 599$. As a survey of
constructions, the paper is written for readers new to the subject; for a
broader introduction we refer to the monographs
\cite{HoradamHadamard, deLauneyFlannery, SeberryYamada, Seberry:Yamada:2020}.

Section~\ref{sec:optimisation} gives definitions and classical
constructions of Hadamard matrices. Section~\ref{sec:alpoge} describes
the matrices at the twelve orders listed above. Section~\ref{sec:sequences}
surveys the known constructions, describes the new matrices, ours and
those of others, and tabulates the open orders up to $10\,000$.
Section~\ref{sec:problems} lists open problems. Every Hadamard matrix constructed 
in this paper is publicly available at \url{https://www.math74.com/Hadamard/} and 
was verified by an exact computation of $HH^{\mathsf T}$.

\section{Hadamard matrices: a rapid overview}\label{sec:optimisation}

The identity and all-ones matrices of order
$n$ are $I_n$ and $J_n$; the all-ones row vector is $j$ and $R$ is the matrix with ones on the back-diagonal, with sizes clear from context; in displayed arrays, $+$ and $-$
stand for $1$ and $-1$.

\begin{theorem}[Hadamard \cite{Hadamard1893}]\label{thm:hadamard-bound}
Let $M = (m_{ij})$ be an $n \times n$ real matrix with $|m_{ij}| \le 1$.
Then $|\det M| \le n^{n/2}$, with equality if and only if
$m_{ij} \in \{\pm 1\}$ for all $i,j$ and $MM\T = nI_n$.
\end{theorem}

The bound follows directly from the inequality
$|\det M| \le \prod_i \|r_i\|$, where $r_1, \dots, r_n$ are the rows of
$M$, together with $\|r_i\| \le \sqrt n$. The inequality was conjectured by
Kelvin and proved by Hadamard; see \cite{MaxDet} for its history. Hadamard drew
attention to real matrices meeting the bound, and constructed examples at orders 12 and 20.

\begin{definition}
A \emph{Hadamard matrix} of order $n$ is an $n \times n$ matrix $H$ with
entries in $\{\pm 1\}$ satisfying $HH\T = nI_n$. We write $H(n)$ for a
Hadamard matrix of order $n$. More generally, a \emph{Butson Hadamard
matrix} $\BH(n,k)$ is an $n \times n$ matrix $H$ with entries in the $k$th
roots of unity and $HH^{\ast} = nI_n$ \cite{Butson}, so that a $\BH(n,2)$ is
a Hadamard matrix.
\end{definition}

Permuting and negating rows or columns preserves the defining
equation. Two Hadamard matrices related by such operations are
\emph{equivalent}, and every Hadamard matrix is equivalent to one whose
first row and column consist of $1$s. Comparing the first three rows of
such a normalised matrix shows that $n \in \{1, 2\}$ or $4 \mid n$.

Sylvester observed that $H(a) \otimes H(b)$ is a Hadamard matrix of
order $ab$, so that $\left(\begin{smallmatrix} + & + \\ + & -
\end{smallmatrix}\right)^{\otimes k}$ is an $H(2^k)$
\cite{Sylvester1867}. The next substantial contribution to the
existence of Hadamard matrices after Hadamard is due to Paley
\cite{Paley1933}. For an odd prime power $q$, let $\chi$ be the
quadratic character of $\F_q$ with the convention that $\chi(0) = 0$, and $Q = (\chi(y - x))_{x,y \in \F_q}$
the \emph{Jacobsthal matrix}. The character sum
$\sum_{z \in \F_q} \chi(z)\chi(z + s) = -1$ for $s \neq 0$ gives
$QQ\T = qI - J$ and $QJ = JQ = 0$, and $Q\T = \chi(-1)\,Q$. A Hadamard
matrix $H$ is \emph{skew} if $H + H\T = 2I$, and a \emph{conference matrix}
of order $n$ is a matrix $C$ with zero diagonal, entries $\pm1$ off the
diagonal, and $CC\T = (n-1)I_n$.

\begin{theorem}[Paley \cite{Paley1933}]\label{thm:paley}
Let
\[
  S = \begin{pmatrix} 0 & j \\ \chi(-1)\, j\T & Q \end{pmatrix}.
\]
If $q \equiv 3 \bmod 4$, then $S$ is skew-symmetric and $I + S$ is
  a skew Hadamard matrix of order $q + 1$. If $q \equiv 1 \bmod 4$, then $S$ is a symmetric conference matrix and
  \[
    S \otimes \begin{pmatrix} + & + \\ + & - \end{pmatrix}
    + I_{q+1} \otimes \begin{pmatrix} + & - \\ - & - \end{pmatrix}
  \]
  is a Hadamard matrix of order $2(q+1)$.
\end{theorem}

The eigenvalues of $Q$ are given by Gauss sums, and the consideration of Gauss and Jacobi
sums gives rise to a large number of \textit{cyclotomic} constructions of Hadamard matrices.
For the theory of cyclotomic numbers on which these depend see Storer or Hall \cite{Storer1967, HallCombinatorialTheory}; for constructions from cyclotomy and character sums see Whiteman
\cite{Whiteman1973}, Spence \cite{Spence1975b} and Yamada
\cite{Yamada1989}; and for the related supplementary difference sets the surveys of
Seberry and Yamada \cite{SeberryYamada} and of Wallis, Street and
Wallis \cite{WallisStreetWallis1972}. Cocyclic and group-developed
constructions are treated in \cite{deLauneyFlannery, HoradamHadamard}.

\subsection{Plug-in constructions}\label{sec:plugin}

Another approach builds $H(4m)$ from four $m \times m$ blocks placed
in a fixed array. Williamson \cite{Williamson1944} took $A, B, C, D$ to
be symmetric circulant $\{\pm1\}$-matrices of order $m$ with
$AA\T + BB\T + CC\T + DD\T = 4mI_m$. Since circulant matrices commute,
the array
\[
  \begin{pmatrix}
     A & B & C & D \\
    -B & A & -D & C \\
    -C & D & A & -B \\
    -D & -C & B & A
  \end{pmatrix}
\]
is then a Hadamard matrix of order $4m$. The Hadamard matrix of order
$92$ was obtained in this way, with $m = 23$
\cite{BaumertGolombHall1962}. Turyn's family supplies Williamson
matrices of order $(q+1)/2$ for every prime power $q \equiv 1 \bmod 4$
\cite{Turyn1972}. Williamson matrices do not exist at every odd order,
however. \DJ okovi\'c showed that there are none of order $35$
\cite{Djokovic35}, and Holzmann, Kharaghani and Tayfeh-Rezaie found
none of orders $47$, $53$ or $59$ \cite{HolzmannKharaghaniTayfehRezaie}.

Goethals and Seidel \cite{GoethalsSeidel1970} removed the symmetry
requirement of the Williamson array by varying the array into which the
circulant matrices are inserted.

\begin{theorem}[Goethals--Seidel \cite{GoethalsSeidel1970}]\label{thm:gs}
Let $A, B, C, D$ be circulant $\{\pm1\}$-matrices of order $m$ with
$AA\T$ $+$ $BB\T$ $+$ $CC\T$ $+$ $DD\T$ $=$ $4mI_m$. Then
\[
  \begin{pmatrix}
     A & BR & CR & DR \\
    -BR & A & D\T R & -C\T R \\
    -CR & -D\T R & A & B\T R \\
    -DR & C\T R & -B\T R & A
  \end{pmatrix}
\]
is a Hadamard matrix of order $4m$.
\end{theorem}

In 2026, Williamson, Yacobi and Zinn-Justin generated Goethals--Seidel
Hadamard matrices of orders up to $252$ by local search,
guided by a transformer trained on the output of that search \cite{WilliamsonYacobiZinnJustin}.

An alternative approach is to allow a constant negative value
off the diagonal in the sum of Gram matrices, which is then absorbed
into a border of $\pm1$ entries. The basic example uses two circulants.

\begin{theorem}[\cite{FletcherGysinSeberry2001}]\label{thm:fgs}
Let $A$ and $B$ be circulant $\{\pm 1\}$-matrices of order $\ell$
satisfying
\begin{equation}\label{eq:lp-gram}
  AA\T + BB\T = (2\ell + 2) I_\ell - 2J_\ell ,
\end{equation}
normalised so that every row of $A$ and of $B$ sums to $1$. Then
\[
 \left(\def\arraystretch{1.2}\begin{array}{cc|cc}
     - & - & j & j\\
     - & + & j & -j\\ \hline
    j\T & j\T & A & B\\
    j\T & -j\T & B\T & -A\T
  \end{array}\right)
\]
is a Hadamard matrix of order $2\ell + 2$.
\end{theorem}

Condition~\eqref{eq:lp-gram} forces the row sums of $A$ and $B$ to be
$\pm 1$, and so $\ell$ is odd. For a $\{\pm1\}$-sequence
$a = (a_0, \dots, a_{m-1})$, the \emph{periodic autocorrelation}
$\PAF_a(s) = \sum_{i=0}^{m-1} a_i a_{i+s}$, with indices mod $m$, is the
inner product of rows $i$ and $i + s$ of the circulant with first row $a$.
The first rows of $A$ and $B$ form a \emph{Legendre pair}: two
$\{\pm1\}$-sequences whose periodic autocorrelations sum to $-2$ at every
non-zero shift. The name
comes from the Legendre sequence of a prime $p \equiv 3 \bmod 4$, whose
periodic autocorrelation is $-1$ at every non-zero shift, so that two
copies of it form a Legendre pair. Legendre pairs are known at infinitely many
lengths and it is conjectured that one exists at every odd length
\cite{DjokovicKotsireas2015}, which would give a Hadamard matrix at
every order $2\ell + 2 \equiv 0 \bmod 4$.

We conclude this section with a construction combining the four circulants of Theorem \ref{thm:gs} with the border and constant negative autocorrelation of Theorem \ref{thm:fgs}.

\begin{theorem}\label{thm:bordered-gs}
Let $A, B, C, D$ be circulant $\{\pm1\}$-matrices of even order $m$ with
\begin{equation}\label{eq:paf-minus-four}
  AA\T + BB\T + CC\T + DD\T = (4m+4)I_m - 4J_m ,
\end{equation}
normalised so that $A$ has row sum $2$ and $B$, $C$, $D$ have row sum
$0$. Then
\[
  \left(\def\arraystretch{1.25}\begin{array}{cccc|cccc}
    + & - & - & - & j & j & j & j \\
    - & + & - & - & j & -j & j & -j \\
    - & - & + & - & j & -j & -j & j \\
    - & - & - & + & j & j & -j & -j \\ \hline
    j\T & j\T & j\T & j\T & A & BR & CR & DR \\
    -j\T & j\T & j\T & -j\T & -BR & A & D\T R & -C\T R \\
    -j\T & -j\T & j\T & j\T & -CR & -D\T R & A & B\T R \\
    -j\T & j\T & -j\T & j\T & -DR & C\T R & -B\T R & A
  \end{array}\right)
\]
is a Hadamard matrix of order $4m + 4$.
\end{theorem}

Multiplying~\eqref{eq:paf-minus-four} by $J_m$ shows that the squares of
the four row sums add to $4$. For even $m$ the row sums are even, so one
of them is $\pm 2$ and the others are $0$. Negating or permuting
$A, B, C, D$ preserves~\eqref{eq:paf-minus-four}.
By analogy with Legendre pairs, the term \emph{Legendre quadruples} seems appropriate 
for four $\{\pm1\}$-sequences of the same length whose periodic
autocorrelations sum to $-4$ at every non-zero shift. The first rows of
$A, B, C, D$ in Theorem~\ref{thm:bordered-gs} form a Legendre quadruple of
even length $m$.

\section{The twelve orders settled in 2026}\label{sec:alpoge}

Each of the twelve orders is of the form $4p$ with $p$ a prime,
$p \equiv 3 \bmod 4$: the primes are $167$, $179$, $223$, $283$, $311$,
$347$, $359$, $419$, $443$, $479$, $487$, $491$. The announcement encodes
the twelve matrices as a single string of signs, together with a short
\texttt{sed} program that assembles them, and does not describe how they
were found \cite{Alpoge2026}. Each assembled matrix is a Goethals--Seidel
array of order $4v$ in four generators $X, Y, Z, W$ over an abelian group
$G$ of order $v$, with a border of width $b$, so that $n = 4v + b$
(Table~\ref{tab:twelve}). 

Over a non-cyclic group, circulants are replaced by group-developed
matrices $[f]_{x,y} = f(y - x)$, for $f \colon G \to \{\pm 1\}$, and $R$
by the matrix of the map $x \mapsto -x$; Theorems~\ref{thm:gs}
and~\ref{thm:bordered-gs} then hold as stated. The periodic autocorrelation
$\PAF_f(d) = \sum_{x \in G} f(x)\, f(x + d)$ is the inner product of rows
$x$ and $x + d$ of $[f]$.
In all twelve matrices
\begin{equation}\label{eq:twelve-paf}
  \sum_{f} \PAF_f(d) =
  \begin{cases}
    -b & \text{if } d \in S \setminus \{0\}, \\
    \phantom{-}4 & \text{if } d \notin S,
  \end{cases}
\end{equation}
for a subgroup $S$ of $G$, where the sum is over the generators $f$ of $X$,
$Y$, $Z$ and $W$. When $b \in \{0, 4\}$, $G$ is cyclic, $S = G$ and
\eqref{eq:twelve-paf} is the hypothesis of Theorem~\ref{thm:gs} or
Theorem~\ref{thm:bordered-gs}; when $b \ge 12$, $S$ is a proper subgroup
of index $b/4 + 1$. Summing~\eqref{eq:twelve-paf} over $d \neq 0$ gives
$\sum_f s_f^2 = 4v + \sum_{d \ne 0} \sum_f \PAF_f(d)$ for the row sums
$s_f$ of the generators.

\begin{table}[ht]
\centering
\caption{The twelve matrices;
$n = 4v + b$ is the matrix order, $b$ the border width, $G$ is a group of order $v$,
the column $\sum \PAF$ gives the values
of~\eqref{eq:twelve-paf} on $S \setminus \{0\}$ and, when $S \ne G$, off
$S$, $U$ is the multiplier group (Definition~\ref{def:multiplier}) and $r$ the number of
its orbits on $G$.}
\label{tab:twelve}
\begin{tabular}{rrlrcrr}
\toprule
$n$ & $b$ & $G$ & $v$ & $\sum \PAF$ & $|U|$ & $r$ \\
\midrule
892  & 0  & $\Z_{223}$ & 223 & $0$ & 3 & 75 \\
1132 & 0  & $\Z_{283}$ & 283 & $0$ & 3 & 95 \\
1244 & 0  & $\Z_{311}$ & 311 & $0$ & 5 & 63 \\
1948 & 0  & $\Z_{487}$ & 487 & $0$ & 9 & 55 \\
1964 & 0  & $\Z_{491}$ & 491 & $0$ & 7 & 71 \\
\midrule
668  & 4  & $\Z_{166}$ & 166 & $-4$ & 1  & 166 \\
716  & 4  & $\Z_{178}$ & 178 & $-4$ & 1  & 178 \\
1676 & 4  & $\Z_{418}$ & 418 & $-4$ & 15 & 42 \\
1772 & 4  & $\Z_{442}$ & 442 & $-4$ & 12 & 50 \\
\midrule
1916 & 12 & $\Z_2^2 \times \Z_{119}$ & 476 & $-12$, $+4$
  & $12 / 4$\textsuperscript{a} & $46 / 128$\textsuperscript{a} \\
1388 & 20 & $\Z_2 \times \Z_3 \times \Z_{57}$ & 342 & $-20$, $+4$ & 6 & 70 \\
1436 & 28 & $\Z_2^2 \times \Z_8 \times \Z_{11}$ & 352 & $-28$, $+4$
  & 5 & 96 \\
\bottomrule
\end{tabular}

\smallskip
\raggedright\footnotesize
\textsuperscript{a}\,$X$ and $W$ are fixed by a group of order $12$,
with $46$ orbits, and $Y$ and $Z$ by a group of order $4$, with $128$
orbits.
\end{table}

\label{sec:multipliers}
\begin{definition}\label{def:multiplier}
Let $U \le \mathrm{Aut}(G)$. A map $f \colon G \to \{\pm1\}$ is
\emph{$U$-invariant}, and $U$ is a \emph{multiplier group} for $f$, if
$f(\sigma x) = f(x)$ for all $\sigma \in U$ and $x \in G$, that is, if $f$
is constant on the $U$-orbits of $G$.
\end{definition}

The five matrices with $b = 0$, at orders $892$, $1132$, $1244$, $1948$ and
$1964$, are given by Theorem~\ref{thm:gs} with circulants of prime order
$m = p$, and their generators share a multiplier group.
In the theory of difference sets a multiplier $t$ of a set $D \subseteq \Z_m$
need only satisfy $tD = D + g_t$ for some $g_t \in \Z_m$ \cite{BJL}; when
$\gcd(|U|, m) = 1$ such a $D$ has a $U$-invariant translate, with the
same periodic autocorrelation.

\begin{lemma}\label{lem:multiplier-paf}
If $a$ is $U$-invariant, then $\PAF_a(ts) = \PAF_a(s)$ for all $t \in U$
and $s \in \Z_m$.
\end{lemma}

If the four generators are $U$-invariant and $U$ has $r$ orbits on
$\Z_m$, each generator is determined by a choice of sign for each orbit. The
multiplier group is a restriction chosen for the search. For $m = p$
prime, $r = 1 + (p-1)/|U|$. At these five orders $|U| = 3, 3, 5, 9, 7$,
and no generator is symmetric.

The four matrices with $b = 4$, at orders $668$, $716$, $1676$ and $1772$, are
given by Theorem~\ref{thm:bordered-gs}, up to the signs of the border, from
Legendre quadruples of length $m = p - 1$. At orders $668$ and $716$ the generators have no multiplier, so all $4m$ signs are free; Best, \DJ okovi\'c, Kharaghani and Ramp observed
that Turyn-type sequences $TT(56)$ and $TT(60)$ would give these two
orders \cite{BestDjokovicKharaghaniRamp}. At orders $1676$ and $1772$,
$m = 418 = 2 \cdot 11 \cdot 19$ and $m = 442 = 2 \cdot 13 \cdot 17$, with
multiplier groups of orders $15$ and $12$. For squarefree
$m = p_1 \cdots p_s$, $r$ is the sum of
$\prod_{i \in I} (p_i - 1)/|U_I|$ over the subsets
$I \subseteq \{1, \dots, s\}$, where $U_I$ is the image of $U$ in
$\prod_{i \in I} \F_{p_i}^{\times}$.

The three matrices with $b \ge 12$ are over non-cyclic groups: at order $1916$,
$G = \Z_2^2 \times \Z_{119}$ and $S = \Z_{119}$; at order $1388$,
$G = \Z_2 \times \Z_3 \times \Z_{57}$ and $S = \Z_{57}$; at order $1436$,
$G = \Z_2^2 \times \Z_8 \times \Z_{11}$ and $S = 2\Z_8 \times \Z_{11}$.
Each is given by the following theorem, with $k = 4$, $6$ and $8$.

\begin{theorem}\label{thm:coset-bordered}
Let $G$ be an abelian group of order $v$, let $S$ be a proper subgroup of $G$
of index $k$, and let $b = 4(k-1)$. Let $\Pi_S$ be the $v \times k$ incidence
matrix of elements and cosets of $S$, and $\Pi = I_4 \otimes \Pi_S$. Let
$X, Y, Z, W$ be group-developed $\{\pm1\}$-matrices over $G$ with
\[
  XX\T + YY\T + ZZ\T + WW\T = (4v + b)I_v - 4k\,\Pi_S\Pi_S\T + 4J_v ,
\]
and let $M$ be the array of Theorem~\ref{thm:gs} in $X, Y, Z, W$.
Let $K$ be a Hadamard matrix of order $b$, and let
$E$ and $F$ be $\{\pm1\}$-matrices of sizes $b \times 4k$ and $4k \times b$ with
\[
  EE\T = 4kI_b, \qquad
  FF\T = 4kI_{4k} - 4(I_4 \otimes J_k), \qquad
  M\,\Pi E\T + \Pi F K\T = 0 .
\]
Then
\[
  H = \begin{pmatrix} K & E\,\Pi\T \\ \Pi F & M \end{pmatrix}
\]
is a Hadamard matrix of order $4v + b$.
\end{theorem}

Since $\Pi_S\Pi_S\T$ is the $\{0,1\}$-matrix of pairs lying in a common coset
of $S$, the condition on $X, Y, Z, W$ is~\eqref{eq:twelve-paf}. The theorem
follows on computing $HH\T$ blockwise, using $\Pi\T\Pi = (v/k)I_{4k}$ and
$\Pi_S J_k \Pi_S\T = J_v$.

At order $1388$ the four generators
are fixed by a common group of order $6$. At order $1436$ they are fixed by
multiplication by the quadratic residues modulo $11$ on the factor
$\Z_{11}$. The sets of entries $-1$ of the four generators have sizes
$k_1, \dots, k_4$. When $b \le 4$ they form \emph{supplementary difference
sets} $4\text{-}\{v; k_1, k_2, k_3, k_4; \lambda\}$: every non-zero element of
$G$ arises exactly $\lambda$ times as a difference of two elements of one of
the sets. When $b \ge 12$ they form a \emph{divisible difference family}
relative to $S$ \cite{MomiharaYamada}, in which each element of
$S \setminus \{0\}$ arises $\lambda_1$ times and each element of $G \setminus S$
arises $\lambda_2 = \lambda_1 + b/4 + 1$ times. In every case
$k_1 + k_2 + k_3 + k_4 - \lambda_1 = v + b/4$, with $\lambda_1 = \lambda$
when $b \le 4$ (Table~\ref{tab:twelve-data}).

For a character $\chi$ of $G$ write $\hat f(\chi) = \sum_{x \in G} f(x)\chi(x)$.
By~\eqref{eq:twelve-paf}, $\sum_f |\hat f(\chi)|^2 = n$ for every character
$\chi$ not trivial on $S$, so a generator has a flat spectrum when
$|\hat f(\chi)|^2 = n/4$ at all such $\chi$. We measure flatness by the mean
$\Phi_f$ of $(4|\hat f(\chi)|^2/n - 1)^2$ over these characters. At $b = 0$,
Parseval's identity gives $v(v-1)\Phi_f = \sum_{d \ne 0} \PAF_f(d)^2 -
(s_f^2 - v)^2/v$, so $\Phi_f$ is essentially the mean square of
$\PAF_f(d)/\sqrt{v}$.

A random sequence has $\Phi_f \approx 1$. To model a typical solution, we
treat the values $|\hat f(\chi)|^2$ at characters in distinct orbits of
$\langle U, -1 \rangle$ as independent exponential variables, as they
approximately are for random sequences, conditioned on the four generators
summing to $n$ at each character. A typical solution then has
$\Phi_f \approx 0.6$, and a typical solution with $|\hat f(\chi)|^2 < n/2$
at every such character has $\Phi_f \approx 0.3$.

The generators of the matrices with $b \ge 4$ are typical, with
$0.27 < \Phi_f < 0.93$ (Table~\ref{tab:twelve-data}). Those with $b = 0$
are much flatter: $|\hat f(\chi)|^2 < n/2$ at every non-trivial character,
and $0.08 < \Phi_f < 0.14$, well below the value $0.3$ that this cutoff
alone gives. By a Chernoff bound under the model, at most $2^{-57}$ (at
$1948$) to $2^{-141}$ (at $1132$) of the solutions with the same multiplier
group are as flat. This suggests that the search at $b = 0$ was restricted
to sequences with nearly flat spectra. At $668$ and $716$, by contrast, the
generators have no multiplier group and are not unusually flat, so there
is no evidence that either restriction was used at these orders.

\begin{table}[ht]
\centering\small
\caption{Row sums of the generators, parameters of their sets of entries
$-1$, and the flatness $\Phi_f$ of each generator.}
\label{tab:twelve-data}
\begin{tabular}{rlll}
\toprule
$n$ & $(s_X, s_Y, s_Z, s_W)$ & $4\text{-}\{v; k_1, k_2, k_3, k_4; \lambda\}$
  & $(\Phi_X, \Phi_Y, \Phi_Z, \Phi_W)$ \\
\midrule
892  & $(11, 11, 11, 23)$ & $\{223; 106, 106, 106, 100; 195\}$ & $(0.09, 0.09, 0.09, 0.09)$ \\
1132 & $(19, 19, 19, 7)$  & $\{283; 132, 132, 132, 138; 251\}$ & $(0.08, 0.09, 0.08, 0.08)$ \\
1244 & $(21, 19, 21, 1)$  & $\{311; 145, 146, 145, 155; 280\}$ & $(0.10, 0.09, 0.10, 0.10)$ \\
1948 & $(17, 1, 17, 37)$  & $\{487; 235, 243, 235, 225; 451\}$ & $(0.14, 0.14, 0.12, 0.12)$ \\
1964 & $(29, 27, 15, 13)$ & $\{491; 231, 232, 238, 239; 449\}$ & $(0.13, 0.11, 0.13, 0.10)$ \\
\midrule
668  & $(2, 0, 0, 0)$ & $\{166; 82, 83, 83, 83; 164\}$  & $(0.36, 0.56, 0.56, 0.44)$ \\
716  & $(2, 0, 0, 0)$ & $\{178; 88, 89, 89, 89; 176\}$  & $(0.50, 0.55, 0.48, 0.59)$ \\
1676 & $(2, 0, 0, 0)$ & $\{418; 208, 209, 209, 209; 416\}$ & $(0.48, 0.92, 0.28, 0.36)$ \\
1772 & $(2, 0, 0, 0)$ & $\{442; 220, 221, 221, 221; 440\}$ & $(0.55, 0.62, 0.41, 0.69)$ \\
\midrule
1916 & $(6, 6, 38, 20)$    & $\{476; 235, 235, 219, 228; 438, 442\}$ & $(0.58, 0.38, 0.34, 0.69)$ \\
1388 & $(6, 10, 26, 24)$   & $\{342; 168, 166, 158, 159; 304, 310\}$ & $(0.49, 0.73, 0.53, 0.45)$ \\
1436 & $(-14, 18, -30, 4)$ & $\{352; 183, 167, 191, 174; 356, 364\}$ & $(0.67, 0.58, 0.54, 0.78)$ \\
\bottomrule
\end{tabular}
\end{table}

\section{Construction of Hadamard matrices of order $\le 10\,000$}\label{sec:sequences}

Theorems~\ref{thm:gs}, \ref{thm:fgs} and~\ref{thm:bordered-gs} reduce
the construction of a Hadamard matrix to the construction of a few
circulant matrices. A circulant is determined by its first row, so the
unknowns are sequences, and the hypotheses of the three theorems say that
the periodic autocorrelations of the first rows sum to $0$, $-2$ and $-4$
respectively at every shift $s \neq 0$. The \emph{non-periodic
autocorrelation} of $a = (a_0, \dots, a_{m-1})$ is
$\NAF_a(s) = \sum_{i=0}^{m-1-s} a_i a_{i+s}$.

Sequences whose \emph{non-periodic} autocorrelations cancel are more
rigid, but they can be concatenated and interleaved to give longer
sequences of the same kind, and most of the recursive constructions rest on
them. The families in use are the following.
\begin{itemize}[nosep]
  \item A \emph{Golay pair} of length $g$ is a pair $a, b$ of
    $\{\pm1\}$-sequences with $\NAF_a(s) + \NAF_b(s) = 0$ for all
    $s \ge 1$ \cite{Golay1961}. A \emph{complex Golay pair} is defined
    in the same way, with entries in $\{\pm1, \pm i\}$ and
    $\NAF_a(s) = \sum_i a_i \overline{a_{i+s}}$.
  \item \emph{Base sequences} $BS(m, n)$ are four $\{\pm1\}$-sequences of
    lengths $m, m, n, n$ whose non-periodic autocorrelations sum to $0$
    at every $s \ge 1$.
  \item \emph{Turyn-type sequences} $TT(n)$ are four $\{\pm1\}$-sequences
    $A, B, C, D$ of lengths $n, n, n, n-1$ with
    $\NAF_A(s)$ $+$ $\NAF_B(s)$ $+$ $2\NAF_C(s)$ $+$ $2\NAF_D(s)$ $=$ $0$ at every $s \ge 1$. The
    concatenations $(C; D)$ and $(C; -D)$, together with $A$ and $B$,
    form base sequences $BS(2n-1, n)$ \cite{Turyn1974}.
  \item \emph{$T$-sequences} of length $t$ are four $\{0, \pm1\}$-sequences,
    exactly one of which is non-zero in each position, whose
    non-periodic autocorrelations sum to $0$ at every $s \ge 1$. Base
    sequences $BS(m, n)$ give $T$-sequences of length $m + n$
    \cite{CooperWallis1972, Turyn1974}.
  \item \emph{$T$-matrices} of order $t$ are four circulant
    $\{0, \pm1\}$-matrices $T_1, \dots, T_4$ of order $t$ with disjoint
    supports and $\sum_i T_iT_i\T = tI_t$. The sums of the periodic autocorrelations of the rows generating the circulants are zero at every $s \geq 1$. $T$-sequences generate $T$-matrices but the generating rows of $T$-matrices are not necessarily $T$-sequences.
  \item \emph{Williamson-type matrices} $\mathrm{WT}(w)$ are four
    $\{\pm1\}$-matrices $A, B, C, D$ of order $w$ with
    $XY\T = YX\T$ for all $X, Y \in \{A, B, C, D\}$ and
    $AA\T + BB\T + CC\T + DD\T = 4wI_w$. Plugged into the Williamson
    array they give an $H(4w)$. The Williamson matrices of
    Section~\ref{sec:plugin} are the circulant symmetric case.
\end{itemize}

\begin{table}[ht]\renewcommand{\arraystretch}{1.2}
\centering\small
\caption{Ingredients of the recursive constructions of Hadamard matrices.}
\label{tab:ingredients}
\begin{tabular}{p{3.2cm}p{2.6cm}p{6.6cm}p{1.6cm}}
\toprule
Ingredient & Used in & Known to exist & Smallest open \\
\midrule
Golay pairs, length $g$
  & $T$-sequences
  & $g = 2^a 10^b 26^c$ \cite{Golay1961, Turyn1974}; no other
    $g \le 105$ \cite{BorweinFerguson2003,EliahouKervaireSaffari1990}.
  & 106 \\
Turyn-type sequences $TT(n)$
  & $T$-sequences
  & even $n \le 46$
    \cite{BestDjokovicKharaghaniRamp, TT_London_Kotsireas_2025,
    LondonKotsireasTT46}.
  & 48 \\
Base sequences $BS(n+1, n)$
  & $T$-sequences
  & $n \le 43$ \cite{Djokovic2010, WangZhu2025}, and
    $n = 2^a 10^b 26^c$ \cite{Djokovic2010}.
  & 44 \\
$T$-sequences, length $t$
  & Cooper--Wallis
  & $g$, $2g + 1$ and $g + g'$ for Golay lengths $g, g'$;
    $2n + 1$ from $BS(n+1, n)$; $3n - 1$ from $TT(n)$; Yang
    multiplication \cite{Yang1989}. All $t \le 100$ except possibly
    $t = 97$ \cite{Djokovic2010}.
  & 97 \\
$T$-matrices, order $t$
  & Cooper--Wallis
  & as for $T$-sequences; $t = 97$, $109$, $197$ \cite{Kumar2026};
    $103$ \cite{Schneider2026}; $151$ \cite{Schumer2026}; $163$ and $181$
    \cite{Jones2026}; $151$ and $163$ also, independently, in this paper. 
  & 127 \\
Williamson-type matrices $\mathrm{WT}(w)$
  & Williamson array, Cooper--Wallis, Seberry--Yamada
  & $(q+1)/2$, $q \equiv 1 \bmod 4$ \cite{Turyn1972};
    $q$ from a symmetric $\mathrm{BH}((q-1)/2, 4)$, $q \equiv 1
    \bmod 4$ \cite{Miyamoto1991}; $2uv$ from $\mathrm{WT}(u)$,
    $\mathrm{WT}(v)$ \cite{SeberryYamada}; the orders below $2000$ in
    \cite[Table~A.1]{SeberryYamada} except $543$ and those resting on
    its Corollary~8.8(1); every odd $w \le 69$, and $73$, $79$, $91$,
    $93$, $95$, $97$, $103$, $109$, $127$, $133$, $181$
    \cite{KharaghaniMohammadianTayfehRezaie2026}; $111$ \cite{Kumar2026};
    $127$, from good matrices \cite{Djokovic1993good}.
  & 71 \\
Skew Hadamard matrices
  & Yamada
  & $q + 1$, $q \equiv 3 \bmod 4$ \cite{Paley1933}; doubling; Spence
    family \cite{Spence1975b}; $4n$ for odd $n < 250$; $276$ and $292$
    \cite{Djokovic2023skew}
  & 356 \\
Symmetric conference matrices
  & Yamada
  & $q + 1$, $q \equiv 1 \bmod 4$ \cite{Paley1933};
    $q^2(q+2) + 1$, $q \equiv 3 \bmod 4$ and $q + 2$ prime powers
    \cite{Mathon1978}; $5 \cdot 9^{2t+1} + 1$ \cite{SeberryWhiteman1988};
    $66$ \cite{Gritsenko2021}; $86$ \cite{Verma2026}.
  & 118\textsuperscript{a} \\
Complex Golay pairs
  & Craigen--Holzmann--Kharaghani
  & 1, 3, 5, 11, 13 and every even $n \leq 26$ except 14
    \cite{BrightKotsireasHeinleGanesh2021}; $2mn$ for complex Golay
    lengths $m, n$ \cite{Craigen1994}; further products
    \cite{CHK2002}.
  & 29 \\
\bottomrule
\end{tabular}

\smallskip
\raggedright\footnotesize
\textsuperscript{a}\,A symmetric conference matrix of order $118$ was
communicated to A.~E.~Brouwer in October 2026 and is listed as known in his
tables \cite{BrouwerSRGtables}, but it is not publicly available. The next
open order is $146$.
\end{table}

A theorem of Cooper and Wallis combines $T$-matrices of order $t$ and
Williamson-type matrices of order $w$ into a Hadamard matrix of order
$4tw$ \cite{CooperWallis1972}. Taking $w = 1$, $T$-sequences of length
$t$ give an $H(4t)$, so base sequences $BS(m, n)$ give an $H(4(m+n))$,
and $TT(n)$ give an $H(4(3n-1))$. The Hadamard matrix of order
$428 = 4 \cdot 107$ was obtained in this way from a $TT(36)$
\cite{Kharaghani428}. Golay pairs enter through Turyn's and Yang's
composition theorems, which combine them with base sequences or
$T$-sequences to give $T$-sequences of greater length
\cite{Turyn1974, Yang1989}. Craigen, Holzmann and Kharaghani showed that
complex Golay pairs of lengths $x$ and $y$ give a Hadamard matrix of
order $8(x + y)$ \cite{CHK2002}. These constructions can be
described with orthogonal designs \cite{GeramitaSeberry}.

Not every known order comes from a family. Twenty-nine orders up to
$10\,000$ are reached by none of the families in
Table~\ref{tab:constructions}, only by matrices found by computer. Four
of these were obtained by \DJ okovi\'c and coauthors from supplementary
difference sets in the Goethals--Seidel array: $764$ \cite{Djokovic764},
$956$ \cite{Djokovic956}, and $1004$ and $2524$
\cite{DjokovicGolubitskyKotsireas}. Twelve are those of
Section~\ref{sec:alpoge}; five, $3732$, $4492$, $5020$, $5788$ and $5932$,
were found by the present authors; and eight, described below, were found
by Schumer, Cooley and Kumar.
Orders first settled by computer but now
reached by a family are not counted as sporadic: Sawade's order $268$
\cite{Sawade268}, and \DJ okovi\'c's order $652 = 4 \cdot 163$
\cite{Djokovic1992} and our order $8476 = 4 \cdot 13 \cdot 163$, which
the theorem of Cooper and Wallis now gives from $T$-matrices of order
$163$. Cati and Pasechnik have implemented the known
constructions up to $1208$ in SageMath
\cite{HadamardDatabase}.

Table~\ref{tab:constructions} lists the constructions used to cover
orders up to $10\,000$, with the number of orders each reaches.
Table~\ref{tab:ingredients} records where the ingredients of the
recursive constructions are known to exist. The remaining recursive
constructions take Hadamard matrices or Williamson-type matrices as
input, and we refer to \cite{SeberryYamada} for them. 

\begin{table}[ht]
\centering\small
\caption{Constructions of Hadamard matrices and the number of orders
$n \le 10\,000$ each reaches. An order reached by several constructions is counted under each.
$\mathrm{WT}(w)$ denotes Williamson-type matrices of order $w$.}
\label{tab:constructions}
\begin{tabular}{p{4.3cm}p{8.6cm}r}
\toprule
Construction & Orders & Count \\
\midrule
\multicolumn{3}{l}{\emph{Direct constructions}} \\[2pt]
Paley I \cite{Paley1933}
  & $q + 1$, $q \equiv 3 \bmod 4$. & 625 \\
Paley II \cite{Paley1933}
  & $2(q + 1)$, $q \equiv 1 \bmod 4$. & 355 \\
Sylvester \cite{Sylvester1867}
  & $2^a$. & 12 \\
Spence \cite{Spence1975b}
  & $4(1 + q + q^2)$, with $1 + q + q^2$ prime $\equiv 3, 5, 7 \bmod 8$
    or $2q^2 + 2q + 3$ a prime power. & 12 \\
Mathon \cite{Mathon1978}
  & $2(q^2(q+2) + 1)$, $q \equiv 3 \bmod 4$, $q + 2$ a prime power.
  & 3 \\
Seberry--Whiteman \cite{SeberryWhiteman1988}
  & $2(5 \cdot 9^{2t+1} + 1)$. & 2 \\
\midrule
\multicolumn{3}{l}{\emph{Recursive constructions}} \\[2pt]
Cooper--Wallis \cite{CooperWallis1972}
  & $4tw$, given $T$-matrices of order $t > 1$ and $\mathrm{WT}(w)$.
  & 1950 \\
Paley and Miyamoto blocks \cite{SeberryYamada1990, Kumar2026}
  & $4st$, given $T$-matrices of order $t > 1$, with $4s - 1 \equiv 3 \bmod 8$
    a prime power; also $6108$, $7532$, $9708$.
  & 1662 \\
Williamson array \cite{Williamson1944}
  & $4w$, given $\mathrm{WT}(w)$. & 1595 \\
Miyamoto II \cite{Miyamoto1991}
  & $\mathrm{WT}(q)$ for $q \equiv 1 \bmod 4$, given a symmetric
    $\mathrm{BH}((q-1)/2, 4)$; reaches $4tq$.
  & 1491 \\
Agaian \cite{Agaian1985}
  & $8ab$, given $H(4a)$, $H(4b)$. & 1230 \\
Kronecker product
  & $ab$, given $H(a)$, $H(b)$. & 1227 \\
Craigen--Holzmann--Kharaghani \cite{CHK2002}
  & $8(x + y)$, given complex Golay pairs of lengths $x$ and $y$.
  & 1001 \\
Craigen--Seberry--Zhang \cite{CraigenSeberryZhang1992}
  & $16abcd$, given $H(4a)$, $H(4b)$, $H(4c)$, $H(4d)$. & 624 \\
Miyamoto I \cite{Miyamoto1991}
  & $4q$, $q \equiv 1 \bmod 4$, given $H(q - 1)$. & 200 \\
Seberry--Yamada \cite{SeberryYamada}
  & $8q$, given $\mathrm{WT}((q-1)/2)$. & 165 \\
Yamada \cite{Yamada1989}
  & $4(q + 2)$, $q \equiv 5 \bmod 8$, given skew $H((q+3)/2)$. & 78 \\
Miyamoto III \cite[Theorem~7]{Miyamoto1991}
  & $8q$, $q \equiv 3 \bmod 4$, with $q$ and $2q - 3$ prime. & 46 \\
Yamada \cite{Yamada1989}
  & $4(q + 2)$, $q \equiv 1 \bmod 8$, given a symmetric conference
    matrix of order $(q+3)/2$. & 34 \\
\midrule
Sporadic (computer search)
  & $764$ \cite{Djokovic764},
    $956$ \cite{Djokovic956}, $1004$ and $2524$
    \cite{DjokovicGolubitskyKotsireas}; $668$, $716$, $892$, $1132$,
    $1244$, $1388$, $1436$, $1676$, $1772$, $1916$, $1948$, $1964$
    \cite{Alpoge2026}; 3732, 4492, 5020, 5788, 5932 (this paper);
    $2092$, $2292$, $2588$ \cite{Cooley2026}; $3676$ \cite{Schumer2026,
    Cooley2026}; $2148$, $2284$,
    $2636$, $3964$ \cite{Kumar2026}.
  & 29 \\
\bottomrule
\end{tabular}
\end{table}

\begin{remark}
Corollary~2 of \cite{Miyamoto1991}, and Corollary~8.8(1) of
\cite{SeberryYamada}, omit the hypothesis that the input be symmetric,
which the construction listed in Table~\ref{tab:constructions} as Miyamoto~II requires. As a
result, order $2076$ is wrongly claimed in \cite[\S7]{Miyamoto1991},
and order $6108$ in \cite[Proposition~3.3]{Djokovic2010}. We construct
a Hadamard matrix of order $2076$ below, and Kumar has constructed one of
order $6108$ \cite{Kumar2026}.
\end{remark}

The known constructions described above give Hadamard matrices at all
but $133$ orders $n \le 10\,000$. By computer search we found Hadamard
matrices at sixteen of these orders. To construct a matrix at order $2076$, we found by computer a
symmetric $\BH(86, 4)$, built from a doubly bordered array and two
matrices of order $42$ with entries in $\{\pm1, \pm i\}$, one circulant
and one back-circulant. It gives Williamson-type matrices of order
$173$, and so the orders $4 \cdot 173 \cdot t$ for odd $t \le 13$, closing six orders otherwise open.
The matrices at orders $3732$, $4492$, $5932$ and $5020$ are given by
Theorem~\ref{thm:gs}, with generators invariant under multiplier groups
of orders $31$, $33$, $39$ and $50$. Those at orders $5788 = 4(4 \cdot 19^2 + 3)$ and
$8476 = 4(4 \cdot 23^2 + 3)$ are given by
Theorem~\ref{thm:coset-bordered}, with $G = \Z_2^2 \times \F_{q^2}$,
$S = \F_{q^2}$ and $q = 19, 23$. We also found $T$-matrices of orders
$t = 151$ and $163$, invariant under multiplication by the subgroup of
order $5$ (resp.\ $9$) of $(\Z/t\Z)^{\times}$. With Williamson matrices of orders $9$,
$11$ and $13$, the theorem of Cooper and Wallis gives the orders $5436$,
$6644$ and $7852$ from $t = 151$, and order $7172$ from $t = 163$. These
$T$-matrices, and these four orders, had been posted earlier by others; we
found ours independently. Schumer gave $T$-matrices of
order $151$ and Hadamard matrices of orders $5436$, $6644$ and $7852$
\cite{Schumer2026}, and Jones gave $T$-matrices of order $163$ and
Hadamard matrices of orders $7172$ and $8476$ \cite{Jones2026}.

Twenty-two further orders were settled by others while this paper was in
preparation; we checked each of their matrices by computing $HH\T$
exactly. Cooley found matrices of orders $2092$, $2292$ and $2588$ by
Theorem~\ref{thm:bordered-gs}, from Legendre quadruples of lengths $522$,
$572$ and $646$, and of order $3676$ by Theorem~\ref{thm:gs} over
$\Z_{919}$ \cite{Cooley2026}; Schumer had earlier given a skew Hadamard
matrix of order $3676$ by the same theorem \cite{Schumer2026}. Kumar found a matrix of order $2284$ by
Theorem~\ref{thm:gs} over $\Z_{571}$, and matrices of orders $2636$ and
$3964$ by Theorem~\ref{thm:coset-bordered} with $k = 4$, over
$\Z_4^2 \times \Z_{41}$ and $\Z_2^2 \times \Z_{247}$ \cite{Kumar2026}.
Kharaghani, Mohammadian, Tayfeh-Rezaie and Zaitsev independently found a
$\BH(1142, 4)$, which also gives order $2284$ \cite{Zaitsev2026}. At order
$2148$ Kumar used a border of width $12$ over $\Z_3 \times \Z_{178}$, with
$S = \Z_{178}$ of index $3$, where the sum of the periodic
autocorrelations is $-12$ on $S \setminus \{0\}$ and $0$ off $S$.

Kumar constructed thirteen further orders without search, from a product
theorem of Seberry and Yamada \cite[Theorem~4]{SeberryYamada1990}: if a
Hadamard matrix of order $4s$ is a $4 \times 4$ array of $s \times s$
circulants $N_{kj}$, and $T_1, \dots, T_4$ are $T$-matrices of order $t$,
then the matrices $A_j = \sum_k T_k \otimes N_{kj}$ satisfy
$\sum_j A_jA_j\T = 4st\,I$, and the Goethals--Seidel array over
$\Z_t \times \Z_s$ gives a Hadamard matrix of order $4st$. For a prime
power $q \equiv 3 \bmod 8$, Paley's matrix of order $q + 1$ is equivalent
to such an array with $s = (q+1)/4$ \cite{Yamada1991, SchmidtIto1999}. With
$T$-matrices of orders $t = 3$, $5$ and $7$ this gives the ten orders
$3156$, $5260$, $5860$, $6036$, $6356$, $6756$, $7364$, $7660$, $8204$ and
$8620$, each of the form $4st$ with $4s - 1$ prime. Miyamoto's matrices of
order $4q$ \cite{Miyamoto1991}, built from near-Williamson matrices
\cite{KharaghaniMohammadianTayfehRezaie2026} or from a periodic Golay pair
\cite{DjokovicKotsireasPGP2015}, can be arranged in the same way, and give
the orders $6108$, $7532$ and $9708$, with $q = 509$, $269$ and $809$.
Finally, order $9172 = 4 \cdot 2293$ follows from Cooley's matrix of order
$2292$ by Miyamoto~I.

Table~\ref{tab:open} lists the remaining $95$ orders $n \le 10\,000$ at which no
Hadamard matrix is currently known. It was produced by applying the
constructions of Table~\ref{tab:constructions} to all orders currently
known, including the twelve of Section~\ref{sec:alpoge} and those
described above. It updates \cite[Table~A.17]{Seberry:Yamada:2020}, which
gives, for each odd $q$, the least $t$ for which a Hadamard matrix of order
$2^t q$ was known.
The smallest open order is now $2396 = 4 \cdot 599$. Every prime appearing in this table is congruent to $3 \bmod 4$. 

\begin{table}[ht]
\centering\small
\caption{Orders $n \le 10\,000$ at which no Hadamard matrix is known,
grouped by the form of $n$. Each entry is the prime $p$; the second
column gives the number of orders in each row.}
\label{tab:open}
\begin{tabular}{lrp{12.5cm}}
\toprule
$n$ & \# & $p$ \\
\midrule
$4p$ & 69 &
    599, 643, 719, 739, 751, 839, 859, 863, 907, 947,
    971, 1019, 1031, 1039, 1051, 1087, 1103, 1187, 1223, 1259,
    1283, 1291, 1319, 1327, 1367, 1423, 1427, 1439, 1451, 1471,
    1487, 1499, 1543, 1559, 1567, 1571, 1579, 1583, 1619, 1663,
    1667, 1699, 1747, 1783, 1787, 1831, 1847, 1871, 1879, 1907,
    1979, 1987, 2027, 2039, 2063, 2083, 2087, 2099, 2203, 2207,
    2251, 2287, 2339, 2347, 2371, 2383, 2399, 2423, 2459 \\
$4\cdot3\,p$ & 7 &
    223, 419, 479, 491, 571, 631, 643 \\
$4\cdot5\,p$ & 6 &
    191, 223, 359, 419, 443, 463 \\
$4\cdot7\,p$ & 3 &
    167, 283, 347 \\
$4\cdot9\,p$ & 1 &
     191 \\
$4\cdot11\,p$ & 1 & 179 \\
$4\cdot13\,p$ & 2 & 167, 179 \\
$8p$ & 6 &
    599, 719, 739, 839, 1031, 1087 \\
\midrule
 & 95 & \\
\bottomrule
\end{tabular}
\end{table}

\section{Open problems}\label{sec:problems}

We collect open problems on Hadamard matrices and on five related
topics, with recent references for each. The open cases in
Table~\ref{tab:ingredients} form a further series of open problems.

\subsection{Hadamard matrices}

\begin{problem}[Legendre pairs]\label{prob:lp}
Does a Legendre pair exist at every odd length $\ell$?
\end{problem}

\begin{problem}[Kotsireas, 2019]\label{prob:lp-pq2}
Let $p$ and $q$ be odd primes, not necessarily distinct, and let
$\left( \frac{\cdot}{p} \right)$ be the Legendre symbol modulo $p$. Do the sequences
$A_{(p,q)} = \left[ 1, q \cdot \left( \frac{1}{p} \right), \ldots, q \cdot \left( \frac{p-1}{p} \right) \right]$, $B_{(p,q)} = \left[ 1, -q \cdot \left( \frac{1}{p} \right), \ldots, -q \cdot \left( \frac{p-1}{p} \right) \right]$
of length $p$ always $q^2$-uncompress to a Legendre pair of length $pq^2$,
that is, are they the $q^2$-compressions of some putative Legendre pair of length $pq^2$? Here the
$q^2$-compression of $(a_0, \dots, a_{pq^2-1})$ is the sequence of length $p$
with $j$th entry $\sum_{i=0}^{q^2-1} a_{j+ip}$. See \cite{KGGG:2026} for
computational results.
\end{problem}

\begin{problem}[Turyn-type sequences]\label{prob:tt}
Does $TT(n)$ exist for every even $n$? Is there an infinite family of
Turyn-type sequences?
\end{problem}

Paley's matrices
of order $q + 1$ with $q \equiv 3 \bmod 4$ are skew, and Seberry
conjectured that skew Hadamard matrices exist at every order divisible
by $4$; see \cite{KoukouvinosSkew}.

\begin{problem}[Skew Hadamard matrices]\label{prob:skew}
Is there a skew Hadamard matrix at every order $4n$?
\end{problem}

\subsection{The maximal determinant problem}

When $n > 2$ and $4 \nmid n$ there is no Hadamard matrix, and one asks
for
\[
  D_n = \max \{ |\det M| : M \in \{\pm1\}^{n \times n} \}.
\]
For odd $n$, Barba's bound $D_n \le B_n = \sqrt{2n-1}\,(n-1)^{(n-1)/2}$
can be attained only if $2n-1$ is a square \cite{Barba}. For
$n \equiv 2 \bmod 4$, the bound
$D_n \le W_n = (2n-2)(n-2)^{(n-2)/2}$ of Ehlich and Wojtas can be
attained only if $2n-2$ is a sum of two squares \cite{Ehlich, Wojtas}.
For $n \equiv 3 \bmod 4$, Ehlich gave a sharper bound $E_n$, which for
$n \ge 63$ can be attained only if $7 \mid n$ \cite{Ehlich}; Cohn showed
that $E_n^2$ is an integer only when $n = 112t^2 \pm 28t + 7$
\cite{Cohn2000}, and Tamura, using the Hasse--Minkowski theorem, that
the smallest order at which $E_n$ could be attained is at least $511$
\cite{Tamura}. The survey \cite{MaxDet} gives proofs and
the state of the problem in 2021.

At $n \equiv 2 \bmod 4$, Momihara, Suda and Xiang
\cite{MomiharaSudaXiang}, and Li, Momihara and Xiang
\cite{LiMomiharaXiang}, constructed infinite families of matrices with
$|\det M| \ge \frac{n-2}{n-1} W_n$; their orders form a set of density
zero. Butbaia et al.\ found the largest known determinants at
$n = 51, 107, 111, 115, 119$ \cite{Butbaia}, and Nu\~nez Ponasso studied
the problem over roots of unity \cite{NunezPonassoRoots}. At
$n \equiv 3 \bmod 4$ the best ratio known at a single order is
$D_n / E_n \ge 0.8122$, at $n = 119$ \cite{Butbaia}, and the best known
along an infinite family is $0.48$ \cite{MaxDet}.

\begin{problem}[Barba bound]\label{prob:barba}
Is $D_n = B_n$ at every $n$ for which $2n - 1$ is a square? Brouwer
\cite{Brouwer} and Whiteman \cite{Whiteman} showed that it is at
$n = 2q^2 + 2q + 1$ for every odd prime power $q$; the smallest
orders left undecided in \cite{MaxDet} are $85$ and $145$.
\end{problem}

\begin{problem}[Ehlich bound]\label{prob:ehlich}
Is $D_n = E_n$ for any $n \equiv 3 \bmod 4$ with $n > 3$?
\end{problem}

\begin{problem}[Near-optimal families]\label{prob:near-optimal}
Is there an infinite set of orders $n \equiv 3 \bmod 4$ along which
$D_n / E_n \to 1$? Does $D_n / W_n \to 1$ as $n \to \infty$ through all
$n \equiv 2 \bmod 4$? For the second, the best general bound is due to
Brent, Osborn and Smith: $D_n / W_n > 1/\pi$ for large $n$ whenever a
Hadamard matrix of order $n - 2$ exists \cite{BrentOsbornProb}.
\end{problem}

\subsection{Butson Hadamard matrices}

A \emph{complex Hadamard matrix} of order $n$ is an $n \times n$ matrix
$H$ with entries of modulus $1$ and $HH^{\ast} = nI$; the character table
of $\Z_n$ is one. The Butson Hadamard matrices $\BH(n,k)$ of
Section~\ref{sec:optimisation} are the complex Hadamard matrices with
entries in the $k$th roots of unity. The inner product of two rows
is a vanishing sum of roots of unity, and the theorem of Lam and Leung
\cite{LamVanishingRoots} shows that a $\BH(n,p^r)$ with $p$ prime and
$n > 1$ exists only if $p \mid n$. A $\BH(n,4)$ gives a Hadamard matrix
of order $2n$ \cite{Turyn}. No such condition arises for $k = 6$, and
Craigen and Sz\H{o}ll\H{o}si constructed a $\BH(p^2, 6)$ for every prime
$p$ \cite[Theorem~1.4.41]{SzollosiThesis}. Tadej and \.Zyczkowski survey
complex Hadamard matrices \cite{TadejZyczkowski}, and Sz\H{o}ll\H{o}si's
thesis treats the Butson case \cite{SzollosiThesis}. A non-existence criterion using 
algebraic number theory is due to Winterhof \cite{Winterhof}. 

The first $\BH(94,4)$ was found by Sz\H{o}ll\H{o}si, from circulant
$\{\pm1\}$-matrices of order $47$ placed in a modification of a block
construction of Kharaghani and Seberry \cite{Szollosi94}, and the first $\BH(118,4)$ by Kharaghani, Mohammadian,
Tayfeh-Rezaie and Zaitsev, from four circulant matrices of order $59$ of which
three are symmetric \cite{KharaghaniMohammadianTayfehRezaie2026}.
Czifra, Matolcsi and Sz\H{o}ll\H{o}si derived a necessary condition for
the existence of a $\BH(n,k)$ from positive definite functions on finite
groups, and constructed the first $\BH(62,6)$, $\BH(82,6)$ and
$\BH(146,6)$ \cite{CzifraMatolcsiSzollosi}.

\begin{problem}[Fourth roots]\label{prob:fourth}
Does a $\BH(n,4)$ exist at every even order $n$? This is a conjecture of
Seberry \cite[Conjecture~1.4.11]{SzollosiThesis}, and it implies the
Hadamard conjecture.
\end{problem}

\begin{problem}[Sixth roots]\label{prob:sixth}
Does a $\BH(2p,6)$ exist for every prime $p$? This is a conjecture of
Sz\H{o}ll\H{o}si \cite[Conjecture~1.4.47]{SzollosiThesis}. According to
\cite{CzifraMatolcsiSzollosi}, the undecided even orders below $100$ are
$46$, $74$, $86$ and $94$, all of this form.
\end{problem}

\subsection{Modular Hadamard matrices}

Marrero and Butson defined an \emph{$m$-modular Hadamard matrix} of
order $n$ to be an $n \times n$ $\{\pm1\}$-matrix $H$ with
$HH\T \equiv nI \bmod m$ \cite{MarreroButson72, MarreroButson73}. The
entries of $HH\T$ lie in $[-n, n]$, so an $m$-modular Hadamard matrix of
order $n < m$ is a Hadamard matrix. The \emph{$m$-modular Hadamard
conjecture} asserts that an $m$-modular Hadamard matrix exists at every
order divisible by $4$; the Hadamard conjecture is equivalent to this
assertion for infinitely many $m$. It is known for $m = 12$
\cite{MarreroButson72, MarreroButson73}, $m = 5$ \cite{LeeSzollosi2014}
and $m = 32$, and so for their divisors. Eliahou and Kervaire settled
$m = 32$ using modular Golay quadruples, which give modular Hadamard
matrices through the Goethals--Seidel array
\cite{EliahouKervaire2001}. Kuperberg settled $m = 7$ and $m = 11$ at
all but finitely many orders \cite{Kuperberg2016}. Eliahou and Kervaire
survey the subject \cite{EliahouKervaire2005}. Eliahou constructed a
$64$-modular Hadamard matrix of order $668$, and announced $64$-modular
Hadamard matrices at all orders $4\ell$ \cite{Eliahou668, Eliahou64}.

\begin{problem}[$128$-modular Hadamard matrices]\label{prob:mod64}
Does a $128$-modular Hadamard matrix exist at every order divisible
by $4$?
\end{problem}

\begin{problem}[$7$-modular Hadamard matrices]\label{prob:mod7}
Does a $7$-modular Hadamard matrix exist at every order divisible by
$4$?
\end{problem}

\subsection{Asymptotic existence}

The asymptotic existence theorems of Seberry and her successors have
the following form: there
are constants $a, b \ge 0$ such that for every odd $q$ there is a
Hadamard matrix of order $2^t q$ for every $t \ge a + b \log_2 q$. The
Hadamard conjecture is the case $a = 2$, $b = 0$. Seberry showed that
one may take $b = 2$ \cite{SeberryAsymptotic}. Craigen improved the bound
using signed groups \cite{CraigenAsymptotic}; $b$ was reduced in stages to
$3/8$ \cite{deLauney2009}, and Livinskyi reduced it further
\cite{Livinskyi}. The introduction of \cite{deLauney2009}
summarises these results. De Launey showed that for every
$\varepsilon > 0$ the odd $q$ for which there is a Hadamard matrix of
order $2^{2 + \lfloor \varepsilon \log_2 q \rfloor} q$ form a set of
positive density \cite{deLauney2009}. Paley's constructions give
$\Omega(x / \log x)$ orders up to $x$ at which a Hadamard matrix exists;
de Launey and Gordon showed that Kronecker products of Paley matrices
give at least
$\frac{x}{\log x} \exp\big((C + o(1))(\log\log\log x)^2\big)$, with
$C = 0.8178\ldots$, and that the other known constructions do not improve
this bound \cite{deLauneyGordon}. Alexeev, Jasper and Mixon showed that
for every $\alpha < 17/92$ and every sufficiently large $n$ there is an
$n \times n$ $\{\pm1\}$-matrix with condition number at most
$1 + n^{-\alpha}$, and that for $4 \nmid n$ the condition number is at
least $1 + c \log n / n$ for an absolute constant $c > 0$
\cite{AlexeevJasperMixon}.

\begin{problem}[Arithmetic progressions]\label{prob:density}
Is there an infinite arithmetic progression $\{a + bk : k \ge 0\}$ of
positive integers such that a Hadamard matrix of order $4(a + bk)$ exists
for every $k$?
\end{problem}

\begin{problem}[Bounded exponent]\label{prob:asymptotics}
Is there an absolute constant $c$ such that a Hadamard matrix of order
$2^t q$ exists for every odd $q$ and every $t \ge c$?
\end{problem}

\subsection{Group-developed matrices and Gauss sums}

Let $G$ be an abelian group of order $n$ and exponent $e$, and let
$f \colon G \to \{\pm1\}$. The matrix $[f(y - x)]_{x, y \in G}$ is a
Hadamard matrix if and only if $|\hat f(\chi)|^2 = n$ for every character
$\chi$ of $G$, where $\hat f(\chi) = \sum_{g \in G} f(g)\chi(g)$ lies in
$\Z[\zeta_e]$. The existence of group-developed Hadamard matrices is thus
a question about algebraic integers of absolute value $\sqrt n$ in a
cyclotomic field, and the same holds for Butson Hadamard
matrices. Gauss sums are the basic source of such integers. The quadratic
Gauss sum underlies Paley's construction, and evaluations of Gauss and
Jacobi sums of higher order give difference sets and strongly regular
graphs \cite{XiangGaussSums}, skew Hadamard difference sets
\cite{FengXiang2012}, the matrices of large determinant of
\cite{MomiharaSudaXiang}, and the group-invariant $\BH(q^2, 6)$ of Ong
and Schmidt \cite{OngSchmidt}. In the other direction, the prime ideal
factorisation of an element of absolute value $\sqrt n$ is constrained,
and this gives the main nonexistence results, by Turyn's self-conjugacy
argument \cite{Turyn1965} and Schmidt's field descent
\cite{Schmidt1999}, developed further by Leung and Schmidt
\cite{SchmidtCirculant}. Ryser conjectured that there is no
circulant Hadamard matrix of order greater than $4$ \cite{Ryser}.
Schmidt's monograph \cite{SchmidtLNM2002} develops these methods in detail.
In October 2026 OpenAI released a manuscript, produced by one of its
models, that claims a proof of Ryser's conjecture; its main theorem is formalised in Lean
\cite{OpenAICirculant2026}.

A difference set $D$ in $G$ is \emph{skew Hadamard} if $G$ is the
disjoint union of $D$, $-D$ and $\{0\}$. The quadratic residues of
$\F_q$, $q \equiv 3 \bmod 4$, form one, and every known example in an
abelian group lies in an elementary abelian group.

\begin{problem}[Skew Hadamard difference sets]\label{prob:skew-ds}
If an abelian group contains a skew Hadamard difference set, must it be
elementary abelian? Momihara showed that one counterexample would give
infinitely many \cite{Momihara2020}.
\end{problem}

\begin{problem}[Momihara, Suda and Xiang]\label{prob:msx}
Momihara, Suda and Xiang constructed three families of matrices $M_n$ of
order $n = q + 1$, with $q \equiv 1 \bmod 4$ a prime power, whose
determinants they expressed through Gauss and Jacobi sums without
evaluating them \cite[Theorem~1.6, Classes (4)--(6)]{MomiharaSudaXiang}.
Is $|\det M_n| \ge 0.8\, W_n$ for every such $n \ge 30$ in each family?
They checked that it is for $n \le 1000$.
\end{problem}

A Hadamard matrix is \emph{cocyclic} over a group $G$ of order $n$ if it
is equivalent to $[\psi(g,h)]_{g,h \in G}$ for a $2$-cocycle
$\psi \colon G \times G \to \{\pm1\}$; when $\psi$ is a coboundary
the matrix is equivalent to a group-developed one. The Sylvester and Paley matrices are cocyclic
\cite{HoradamHadamard, deLauneyStaffordPaleyI}.

\begin{problem}[Cocyclic Hadamard matrices]\label{prob:cocyclic}
Is there a cocyclic Hadamard matrix at every order divisible by $4$?
This is a conjecture of de Launey and Horadam \cite{HoradamHadamard}.
\end{problem}

\section*{Statement of AI Use}

The authors used Anthropic's Claude (Claude Opus 5.5, via Claude Code,
September--October 2026) extensively in this work. Claude wrote the Python
and C++ programs used to analyse the twelve matrices of~\cite{Alpoge2026},
to check constructions claimed in the literature and the matrices found by
others described in Section~\ref{sec:sequences}, and to construct the
symmetric $\BH(86,4)$, the $T$-matrices of orders $151$ and $163$, and the
Hadamard matrices of orders $5788$ and $8476$ of
Section~\ref{sec:sequences}. Theorem~\ref{thm:coset-bordered} was formulated
from this analysis, then checked and rewritten by hand. Sections of the paper were drafted using
Claude, then checked and rewritten by hand. Each section was
proof-read by all four authors. Every Hadamard matrix constructed in this paper is publicly available at \url{https://www.math74.com/Hadamard/} and was verified by an exact computation of $HH^{\mathsf T}$. The authors take full responsibility for all claims made in this paper.

\bibliographystyle{abbrv}
\bibliography{Biblio2020}

\end{document}